# Cumulative distribution function estimation under interval censoring case 1

### Elodie Brunel*

*University Montpellier 2*
*I3M UMR CNRS 5149, cc. 051, place E. Bataillon*
*34095 Montpellier cedex, France*
*e-mail:* ebrunel@math.univ-montp2.fr

### Fabienne Comte

*University Paris Descartes*
*MAP5 UMR CNRS 8145, 45 rue des Saints Pères*
*75270 Paris cedex 6, France*
*e-mail:* fabienne.comte@parisdescartes.fr

**Abstract:** We consider projection methods for the estimation of the cumulative distribution function under interval censoring, case 1. Such censored data also known as current status data, arise when the only information available on the variable of interest is whether it is greater or less than an observed random time. Two types of adaptive estimators are investigated. The first one is a two-step estimator built as a quotient estimator. The second estimator results from a mean square regression contrast. Both estimators are proved to achieve automatically the standard optimal rate associated with the unknown regularity of the function, but with some restriction for the quotient estimator. Simulation experiments are presented to illustrate and compare the methods.

## Contents

*Corresponding author.

## 1. Introduction

Let $X$ be a survival time with unknown cumulative distribution function (cdf) $F$. In the interval censoring case 1 model, we are not able to observe the survival time $X$. Instead, an observation consists of the pair $(U, \delta)$ where $U$ is an examination time and $\delta$ is the indicator function of the event $(X \leq U)$. Roughly speaking, the only knowledge about the variable of interest $X$ is wether it has occurred before $U$ or not. Early examples of such interval censoring can be found in demography studies, see Diamond and McDonald (1991). In epidemiology, these censoring schemes also arise for instance in AIDS studies or more generally in the study of infectious diseases when the infection time is an unobservable event. We assume that $U$ is independent of $X$, that $F$ has density $f$ and that the cdf G of $U$ has density $g$. Such data, also known as current status data, may remind us right-censored data where the observed data is the pair $(\min(X, C), \mathbb{I}_{(X \leq C)})$ where $C$ is a censoring variable. However, the estimation procedure in these two censoring models is substantially different. In the right-censoring model, the Kaplan and Meier (1958) estimator is well studied and is asymptotically normal at the rate $\sqrt{n}$. Nevertheless, current status data have been studied by many authors in the last two decades, see Jewell and van der Laan (2004) for a state of the art. In the interval censoring model, the nonparametric maximum likelihood estimator (NPMLE) of the survival function is proved to be uniformly consistent, pointwise convergent to a nonnormal asymptotic distribution at the rate $n^{1/3}$ in Groeneboom and Wellner (1992). In van de Geer (1993), it is also established that the NPMLE converges at rate $n^{-1/3}$ in $\mathbb{L}^2$-norm.

Recent extensions take two directions. First, more general contexts are considered. For example, van der Vaart and van der Laan (2006) build nonparametric estimates of the survival function for current status data in presence of time dependent and high dimensional covariates: they provide limit central theorems with rate $n^{1/3}$ and nonstandard limiting processes. The second direction aims at proposing smooth estimates that may take into account the possible smoothness of the survival function. Indeed, the NPMLE estimator is a piecewise constant function. The locally linear smoother proposed by Yang (2000), contrary to the NPMLE may be non monotone, but it has a better convergence rate than the NPMLE when the density $f$ is smooth and the kernel function and the bandwidth are properly chosen. In the same spirit, Ma and Kosorok (2006) introduce an adaptive modified penalized least square estimator built with smoothing splines but their main objective is the study of semiparametric models. They have in mind the same type of penalization device that we present here, but their penalty functions contain many complicated terms that would be difficult to estimate.

Here, we also pursue the search for smooth (or piecewise smooth) adaptive estimators. We present two different penalized minimum contrast estimators built on trigonometric, polynomial or wavelet spaces whose associated penalty terms are really simple; the minimization of the penalized contrast function allows to choose a space that leads to both a non asymptotic automatic squared bias/variance compromise and to an asymptotic optimal convergence rate according to the regularity of the function $F$ in term of Besov spaces. An interesting feature of the procedure is that the estimators and their study is made straightforward by the most powerful Talagrand (1996) inequality for empirical centered processes. We also use technical properties proved in a regression framework by Baraud *et al.* (2001) and Baraud. (2002) for the mean-square estimator. Globally, the available tools and algorithms for adaptive density and regression estimation make our solution easy to study and to implement.

The plan of the paper is as follows. Section 2 introduces the quotient and the regression estimators, after the description of the lifetimes model. We also give a detailed description of the projection spaces with their main properties. Then, we study one projection estimator of the density of the failure times which have occurred before the examination time in Section 3. Both convergence and adaptation results are given. This estimator is then applied to the estimation of the cumulative distribution function *via* a quotient construction. Section 4 describes a direct adaptive procedure to estimate the distribution function based on a mean square regression contrast. Simulations compare both approaches in Section 5. We use as a benchmark the NPMLE and also the simple piecewise constant estimator proposed by Birgé (1999). Lastly, most proofs and technical lemmas are deferred to Section 6.

## 2. Definition of the estimators

### *2.1. Model and assumptions*

Let $(U_1, \delta_1), \cdots (U_n, \delta_n)$ be a sample of the pair $(U, \delta)$ where $\delta_i = \mathbb{I}_{(X_i \leq U_i)}$ and the sequences $(U_i)_{1 \leq i \leq n}$ and $(X_i)_{1 \leq i \leq n}$ are independent. We are interested in the estimation of the distribution function $F$ of the lifetime $X$ on a compact set $A$ only. By rescaling the data, we take $A = [0, 1]$ without loss of generality.

The compact set is considered as fixed in the theory, even if in practice, it is determined from the data. Remember that we denote by $f$ and $F$ the density and the cumulative distribution function of the unobserved lifetime $X$ and $g$ and $G$ those of the examination time $U$. A function of interest is the density $\psi$ of the $U_i$ restricted to the individuals for which $\delta_i = 1$ defined by:

$$\psi(x) = F(x)\,g(x) \tag{2.1}$$

It is clear that this equation provides a way to build an estimator of $F$. This approach is developed in Section 3. The censoring mechanism is such that the conditional law of $\delta = \mathbb{I}_{(X \leq U)}$ given $U = u$ is a Bernoulli law with parameter $F(u)$ and as a consequence we have:

$$\mathbb{E}(\delta|U = u) = F(u) \tag{2.2}$$

This relation will lead to define a direct mean-square estimator of $F$.

Both strategies require the following assumption:

**[A1]** The density $g$ of the random time $U$ is lower and upper bounded on $A$ so that there exist real finite constants $g_0 > 0$ and $g_1 > 0$ such that for all $x \in A$, $g_0 \leq g(x) \leq g_1$.

### 2.2. Definition of the estimators

Assume that we have at our disposal a collection of finite dimensional spaces of functions, denoted by $(S_m)_{m \in \mathcal{M}_n}$, satisfying the following assumption:

$(\mathcal{H}_1)$ $(S_m)_{m \in \mathcal{M}_n}$ is a collection of finite-dimensional linear sub-spaces of $\mathbb{L}^2([0,1])$, with dimension $\dim(S_m) = D_m$ such that $D_m \leq n$, $\forall m \in \mathcal{M}_n$ and satisfying:
$$\exists \Phi_0 > 0, \forall m \in \mathcal{M}_n, \forall t \in S_m, \|t\|_\infty \leq \Phi_0 \sqrt{D_m} \|t\|. \tag{2.3}$$
where $\|t\|^2 = \int_0^1 t^2(x)dx$, for $t$ in $\mathbb{L}^2([0,1])$.

#### 2.2.1. Quotient estimator

As already mentioned, the first strategy requires to estimate $\psi$ and $g$. The estimator $\tilde{g}$ of $g$ is chosen as the adaptive density estimator defined in Massart (2007), Chapter 7, namely: $\tilde{g} = \hat{g}_{\hat{m}_g}$ where $\hat{g}_m = \arg\min_{t \in S_m} \gamma_n^g(t)$,

$$\gamma_n^g(t) = \|t\|^2 - \frac{2}{n} \sum_{i=1}^n t(U_i),$$

and

$$\hat{m}_g = \arg\min_{m \in \mathcal{M}_n} \gamma_n^g(\hat{g}_m) + \mathrm{pen}_g(m). \tag{2.4}$$

with $\mathrm{pen}_g(m) = \kappa \Phi_0^2 D_m / n$.

For the estimation of $\psi$, we consider the following contrast function

$$\gamma_n^\psi(t) = \|t\|^2 - \frac{2}{n}\sum_{i=1}^n \delta_i t(U_i). \tag{2.5}$$

Let then
$$\hat{\psi}_m = \arg\min_{t \in S_m} \gamma_n^\psi(t). \tag{2.6}$$

Then we define $\tilde{\psi} = \hat{\psi}_{\hat{m}}$ where

$$\hat{m} = \arg\min_{m \in \mathcal{M}_n} [\gamma_n^\psi(\hat{\psi}_m) + \text{pen}^\psi(m)].$$

The penalty function will be motivated and defined later. The contrasts $\gamma_n^g$ and $\gamma_n^\psi$ are both found as empirical versions of the $\mathbb{L}^2$ distance between a function $t$ in $S_m$ and the function of interest ($g$ or $\psi$). To see this, take the expectation of e.g. $\gamma_n^\psi$:

$$\mathbb{E}(\gamma_n^\psi(t)) = \|t\|^2 - 2\langle t, \psi\rangle = \|t - \psi\|^2 - \|\psi\|^2$$

with $\langle t, s\rangle = \int t(x)s(x)dx$. This illustrates that minimizing $\gamma_n^\psi$ is likely to provide a function $t$ that minimizes in mean $\|t - \psi\|^2$ and thus estimate $\psi$, on the space $S_m$.

Now, the adaptive estimators $\tilde{\psi}$ of $\psi$ and $\tilde{g}$ of $g$ are defined, and we can use Equality (2.1) to build a quotient estimator of the distribution function $F$ by setting

$$\tilde{F}(x) = \begin{cases} 0 & \text{if } \tilde{\psi}(x)/\tilde{g}(x) < 0 \\ \dfrac{\tilde{\psi}(x)}{\tilde{g}(x)} & \text{if } 0 \leq \tilde{\psi}(x)/\tilde{g}(x) \leq 1 \\ 1 & \text{if } \tilde{\psi}(x)/\tilde{g}(x) > 1 \end{cases} \tag{2.7}$$

*2.2.2. Regression estimator*

On the other hand, a direct estimator of the cdf $F$ can be obtained by considering the following mean-square contrast:

$$\gamma_n^{\text{MS}}(t) = \frac{1}{n}\sum_{i=1}^n [\delta_i - t(U_i)]^2 \tag{2.8}$$

In this case, we set
$$\hat{F}_m = \arg\min_{t \in S_m} \gamma_n^{MS}(t) \tag{2.9}$$

in the sense that we always can compute a vector $(\hat{F}_m(U_1), \ldots, \hat{F}_m(U_n))$ as the orthogonal projection of the vector $(\delta_1, \ldots, \delta_n)$ on the sub-space of $\mathbb{R}^n$ defined by $\{(t(U_1), \ldots, t(U_n)), t \in S_m\}$. Then we define $\hat{F}_{\hat{m}_0}$ by:

$$\hat{m}_0 = \arg\min_{m \in \mathcal{M}_n} \{\gamma_n^{MS}(\hat{F}_m) + \text{pen}^{MS}(m)\}, \tag{2.10}$$

with
$$\text{pen}^{MS}(m) = \kappa_0 \frac{D_m}{n}. \tag{2.11}$$

where $\kappa_0$ is a numerical constant.

### *2.2.3. Remark about the NPMLE estimator*

In the above setting, it is conceivable to consider the log-likelihood contrast $\gamma_n^{MLE}(t) = (1/n) \sum_{i=1}^{n} (\delta_i \log(t(U_i)) + (1-\delta_i) \log(1-t(U_i)))$. If $t$ is supposed to be a piecewise constant function with jumps only at the observed points, then the NPMLE $\hat{F}_n$ which maximizes $\gamma_n^{MLE}$ can be obtained by the *max-min* formula, see Jewell and van der Laan (2004):

$$\hat{F}_n(U_{(i)}) = \max_{j \leq i} \min_{k \geq i} \frac{\sum_{m=j}^{k} \delta_{(m)}}{k-j+1} \tag{2.12}$$

Most results are essentially of asymptotic nature for the NPMLE as already mentioned. Nevertheless, it is of interest for setting the benchmark to include it in our simulation study, see Section 5. The advantage is that no adaptation is required, but the NPMLE is a piecewise constant function. Now, another approach is to consider the histogram-type estimator introduced by Birgé (1999). Let $(I_j)_{1 \leq j \leq D} = ([a_{j-1}, a_j[)_{1 \leq j \leq D}$ be a partition of $[0,1]$ and let us consider a piecewise constant function $t = \sum_{j=1}^{D} a_j \mathbb{I}_{I_j}$. If one looks for such a function that maximizes the contrast $\gamma_n^{MLE}$, one finds the estimator given by Birgé (1999): $\hat{F}_D = \sum_{j=1}^{D} \hat{a}_j^{MLE} \mathbb{I}_{I_j}$ with

$$\hat{a}_j^{MLE} = \frac{1}{N_j} \sum_{i=1}^{n} \delta_i \mathbb{I}_{I_j}(U_i), \quad \text{if} \ \ N_j = \sum_{i=1}^{n} \mathbb{I}_{I_j}(U_i) \neq 0,$$

and $\hat{a}_j^{MLE} = 0$ otherwise for $j = 1, \ldots, D$. If $D$ is of order $n^{1/3}$, Birgé (1999) gives weak assumptions ensuring that the $\mathbb{L}^1$-risk $\mathbb{E}(\int_0^1 |\hat{F}_D(x) - F(x)| dx)$ is of order $n^{-1/3}$. A thinner model selection strategy may be developed to take a possible higher regularity of $F$ into account.

But the contrasts proposed here have the advantage that the empirical processes to be controlled are linear with respect to the functions $t$, a property that the NPMLE estimator would not share. This would make the theoretical study more technical, and the estimation algorithms difficult to implement for general bases. Moreover, Hellinger-type risk would have to be considered, in the same way as in Birgé and Rozenholc (2006). This explains why we rather consider the contrasts $\gamma_n^\psi$ and $\gamma_n^{MS}$.

Before studying both estimators, let us give some examples of collections $(S_m)_{m \in \mathcal{M}_n}$.

### *2.3. Spaces of approximation*

The main assumption is described by $(\mathcal{H}_1)$. In this setting, an orthonormal basis of $S_m$ is denoted by $(\varphi_\lambda)_{\lambda \in \Lambda_m}$ where $|\Lambda_m| = D_m$. Let us mention that it follows from Birgé and Massart (1997) that Property (2.3) in the context of $(\mathcal{H}_1)$ is

equivalent to

$$\exists \Phi_0 > 0, \|\sum_{\lambda \in \Lambda_m} \varphi_\lambda^2\|_\infty \leq \Phi_0^2 D_m. \tag{2.13}$$

Moreover, for some results we need the following additional assumption:

($\mathcal{H}_2$) $(S_m)_{m \in \mathcal{M}_n}$ is a collection of nested models, we denote by $\mathcal{S}_n$ the space belonging to the collection, such that $\forall m \in \mathcal{M}_n, S_m \subset \mathcal{S}_n$. We denote by $N_n$ the dimension of this nesting space: $\dim(\mathcal{S}_n) = N_n$ ($\forall m \in \mathcal{M}_n, D_m \leq N_n$).

We consider more precisely the following examples:

- [T] *Trigonometric spaces*: $S_m$ is generated by $\{\,1, \sqrt{2}\cos(2\pi j x), \sqrt{2}\sin(2\pi j x)$ for $j = 1, \ldots, m\,\}$, $D_m = 2m+1$ and $\mathcal{M}_n = \{1, \ldots, [n/2]-1\}$.
- [P] *Regular piecewise polynomial spaces*: $S_m$ is generated by $m(r+1)$ polynomials, $r+1$ polynomials of degree $0, 1, \ldots, r$ on each subinterval $[(j-1)/m, j/m]$, for $j = 1, \ldots m$, $D_m = (r+1)m$, $m \in \mathcal{M}_n = \{1, 2, \ldots, [n/(r+1)]\}$. For example, consider the orthogonal collection in $\mathbb{L}^2([-1, 1])$ of Legendre polynomials $Q_k$, where the degree of $Q_k$ is equal to $k$, $|Q_k(x)| \leq 1, \forall x \in [-1, 1]$, $Q_k(1) = 1$ and $\int_{-1}^1 Q_k^2(u)du = 2/(2k+1)$. Then the orthonormal basis is given by $\varphi_{j,k}(x) = \sqrt{m(2k+1)}Q_k(2mx - 2j + 1)\mathbb{I}_{[(j-1)/m, j/m[}(x)$ for $j = 1, \ldots, m$ and $k = 0, \ldots, r$, with $D_m = (r+1)m$. In particular, the histogram basis corresponds to $r = 0$ and is simply defined by $\varphi_j(x) = \sqrt{D_m}\mathbb{I}_{[(j-1)/D_m, j/D_m]}(x)$ and $D_m = m$. We denote by [DP] the collection of piecewise polynomials corresponding to dyadic subdivisions with $m = 2^q$ and $D_m = (r+1)2^q$.
- [W] *Dyadic wavelet generated spaces* with regularity $r$ and compact support, as described e.g. in Donoho and Johnstone (1998).

All those spaces satisfy ($\mathcal{H}_1$), with for instance $\Phi_0 = \sqrt{2}$ for collection [T] and $\Phi_0 = \sqrt{2r+1}$ for collection [P]. Moreover, [T], [DP] and [W] satisfy ($\mathcal{H}_2$) since they are nested with $\mathcal{S}_n$ being the space with the largest dimension in the collection.

## 3. Study of the quotient estimator

Our aim is to estimate the cdf $F$ from the observations $(\delta_i, U_i)$, $i = 1, \ldots, n$.

### 3.1. Convergence results for one estimator

An explicit expression of the estimator follows from definition (2.5)–(2.6) by using the orthonormal basis $(\varphi_\lambda)_{\lambda \in \Lambda_m}$ of $(S_m)$ described in ($\mathcal{H}_1$):

$$\hat{\psi}_m = \sum_{\lambda \in \Lambda_m} \hat{a}_\lambda \varphi_\lambda \text{ with } \hat{a}_\lambda = \frac{1}{n}\sum_{i=1}^n \delta_i \varphi_\lambda(U_i). \tag{3.1}$$

We define also $\psi_m$ as the orthogonal projection of $\psi$ on $S_m$. We can write

$$\psi_m = \sum_{\lambda \in \Lambda_m} a_\lambda \varphi_\lambda \text{ with } a_\lambda = \int_0^1 \varphi_\lambda(x)\psi(x)dx. \tag{3.2}$$

The rate of the estimator $\hat{\psi}_m$ of $\psi$ is quite easy to derive. Indeed, it follows from (3.1), (3.2) and Pythagoras theorem that

$$\begin{aligned}
\|\psi - \hat{\psi}_m\|^2 &= \|\psi - \psi_m\|^2 + \|\psi_m - \hat{\psi}_m\|^2 = \|\psi - \psi_m\|^2 + \sum_{\lambda \in \Lambda_m} (a_\lambda - \hat{a}_\lambda)^2 \\
&= \|\psi - \psi_m\|^2 + \sum_{\lambda \in \Lambda_m} \left( \frac{1}{n} \sum_{i=1}^n \delta_i \varphi_\lambda(U_i) - \int_0^1 \psi(x) \varphi_\lambda(x) dx \right)^2.
\end{aligned}$$

Therefore

$$\begin{aligned}
\mathbb{E}(\|\psi - \hat{\psi}_m\|^2) &= \|\psi - \psi_m\|^2 + \sum_{\lambda \in \Lambda_m} \text{Var}\left( \frac{1}{n} \sum_{i=1}^n \delta_i \varphi_\lambda(U_i) \right) \\
&= \|\psi - \psi_m\|^2 + \frac{1}{n} \sum_{\lambda \in \Lambda_m} \text{Var}\left( \delta_1 \varphi_\lambda(U_1) \right) \\
&\leq \|\psi - \psi_m\|^2 + \frac{1}{n} \mathbb{E}\left[ \left( \sum_{\lambda \in \Lambda_m} \varphi_\lambda^2(U_1) \right) \delta_1 \mathbb{I}_{(U_1 \leq 1)} \right] \\
&\leq \|\psi - \psi_m\|^2 + \frac{\Phi_0^2 D_m}{n} \mathbb{E}(\delta_1 \mathbb{I}_{(U_1 \leq 1)})
\end{aligned}$$

with (2.13). This can be summarized by the following Proposition:

**Proposition 3.1.** *Consider the model described in Section 2.1 and the estimator $\hat{\psi}_m = \arg\min_{t \in S_m} \gamma_n^\psi(t)$ where $\gamma_n^\psi(t)$ is defined by (2.5) and $S_m$ is a $D_m$-dimensional linear space in a collection satisfying $(\mathcal{H}_1)$. Then*

$$\mathbb{E}(\|\psi - \hat{\psi}_m\|^2) \leq \|\psi - \psi_m\|^2 + \frac{\Phi_0^2 D_m}{n} \mathbb{E}(\delta_1 \mathbb{I}_{(U_1 \leq 1)}). \tag{3.3}$$

Inequality (3.3) gives the asymptotic rate for one estimator if we consider that $\psi$ belongs to a Besov space $\mathcal{B}_{\alpha_\psi,p,\infty}([0,1])$ with finite Besov norm denoted by $|\psi|_{\alpha_\psi,p}$. For a precise definition of those notions we refer to DeVore and Lorentz (1993) Chapter 2, Section 7, where it is also proved that $\mathcal{B}_{\alpha_\psi,p,\infty}([0,1]) \subset \mathcal{B}_{\alpha_\psi,2,\infty}([0,1])$ for $p \geq 2$. This justifies that we now restrict our attention to $\mathcal{B}_{\alpha_\psi,2,\infty}([0,1])$.

Then the following (standard) rate is obtained:

**Corollary 3.1.** *Consider the model described in Section 2.1 and the estimator $\hat{\psi}_m = \arg\min_{t \in S_m} \gamma_n^\psi(t)$ where $\gamma_n^\psi(t)$ is defined by (2.5) and $S_m$ is a $D_m$-dimensional linear space in collection* [T], [P], *or* [W]. *Assume moreover that $\psi$* belongs to $\mathcal{B}_{\alpha_\psi,2,\infty}([0,1])$ with $r > \alpha_\psi > 0$ and choose a model with $m = m_n$ such that $D_{m_n} = O(n^{1/(2\alpha_\psi+1)})$, then

$$\mathbb{E}(\|\psi - \hat{\psi}_{m_n}\|^2) = O\left(n^{-\frac{2\alpha_\psi}{2\alpha_\psi+1}}\right). \tag{3.4}$$

**Remark 3.1.** The bound $r$ stands for the regularity of the basis functions for collections [P] and [W]. For the trigonometric collection [T], no upper bound for the unknown regularity $\alpha_\psi$ is required.

*Proof.* The result is a straightforward consequence of the results of DeVore and Lorentz (1993) and of Lemma 12 of Barron *et al.* (1999), which imply that $\|\psi - \psi_m\|$ is of order $D_m^{-\alpha_\psi}$ in the three collections [T], [P] and [W], for any positive $\alpha_\psi$. Thus the minimum order in (3.3) is reached for a model $S_{m_n}$ with $D_{m_n} = O([n^{1/(1+2\alpha_\psi)}])$, which is less than $n$ for $\alpha_\psi > 0$. Then, if $\psi \in \mathcal{B}_{\alpha_\psi,2,\infty}([0,1])$ for some $\alpha_\psi > 0$, we find the standard nonparametric rate of convergence $n^{-2\alpha_\psi/(1+2\alpha_\psi)}$. □

### 3.2. Adaptive estimator of the density $\psi$

The penalized estimator is defined in order to ensure an automatic choice of the dimension. Indeed, it follows from Corollary 3.1 that the optimal dimension depends on the unknown regularity $\alpha_\psi$ of the function to be estimated in the asymptotic setting and more generally on the unknown constants involved in the squared-bias/variance terms. Then we define

$$\hat{m} = \arg\min_{m \in \mathcal{M}_n}[\gamma_n^\psi(\hat{\psi}_m) + \mathrm{pen}^\psi(m)]$$

where the penalty function $\mathrm{pen}^\psi$ is determined in order to lead to the choice of a "good" model. First, we apply some Talagrand (1996) type inequality to the linear empirical process defined by

$$\nu_n(t) := \frac{1}{n}\sum_{i=1}^{n}\left(\delta_i t(U_i) - \langle t, \psi\rangle\right). \tag{3.5}$$

Then, by using the decomposition of the contrast given by

$$\gamma_n^\psi(t) - \gamma_n^\psi(s) = \|t - \psi\|^2 - \|s - \psi\|^2 - 2\nu_n(t-s), \tag{3.6}$$

we easily derive the following result:

**Theorem 3.1.** *Consider the model described in Section 2.1 and the estimator $\hat{\psi}_m = \arg\min_{t \in S_m} \gamma_n^\psi(t)$ where $\gamma_n^\psi(t)$ is defined by (2.5) and $S_m$ is a $D_m$-dimensional linear space in a collection satisfying $(\mathcal{H}_1)$ and $(\mathcal{H}_2)$. Then the estimator $\hat{\psi}_{\hat{m}}$ with $\hat{m}$ defined by*

$$\hat{m} = \arg\min_{m \in \mathcal{M}_n}[\gamma_n^\psi(\hat{\psi}_m) + \mathrm{pen}^\psi(m)]$$

*and*

$$\mathrm{pen}^\psi(m) \geq \kappa \Phi_0^2 \left( \int_0^1 \psi(x) dx \right) \frac{D_m}{n}$$

*where $\kappa$ is a universal constant, satisfies*

$$\mathbb{E}(\|\hat{\psi}_{\hat{m}} - \psi\|^2) \leq \inf_{m \in \mathcal{M}_n} \left( 3\|\psi - \psi_m\|^2 + 4\mathrm{pen}^\psi(m) \right) + \frac{C}{n}, \qquad (3.7)$$

*where $C$ is a constant depending on $\Phi_0$ and on $\int_0^1 \psi(x)dx$.*

As it is clear from Theorem 3.1, only a lower bound for the penalty is provided. As $\mathrm{pen}^\psi(.)$ also appears in the risk bound (3.7), we should not take it much larger than the order $D_m/n$ because then, the $\mathbb{L}^2$ error would increase and the resulting rate would not be the optimal one. On the other hand, no result is available for smaller penalties. This explains in particular why it is possible to keep the asymptotic rate unchanged when increasing the constant $\kappa$ only.

Also, if we choose $\mathrm{pen}^\psi(m) = \kappa \Phi_0^2 \left( \int_0^1 \psi(x)dx \right) (D_m/n)$, it follows from (3.7) that the adaptive estimator automatically makes the squared-bias/variance compromise and from an asymptotic point of view, reaches the optimal rate, provided that the constant in the penalty is known. Note that Inequality (3.7) is nevertheless non-asymptotic.

**Remark 3.2.** In practice, the constant in the penalty, denoted above by $\kappa$, is found by simulation experiments taking into account very different types of functions $\psi$. See examples of such a work in Birgé and Rozenholc (2006) or Comte and Rozenholc (2004).

The penalty given in Theorem 3.1 cannot be used in practice since it depends on the unknown quantity

$$\int_0^1 \psi(x)dx = \mathbb{E}(\delta_1 \mathbb{I}_{(U_1 \leq 1)}).$$

A simple solution is to use that $\int_0^1 \psi(x)dx \leq 1$; it follows that Inequality (3.7) would hold for a penalty defined by $\mathrm{pen}^\psi(m) = \kappa \Phi_0^2 D_m/n$. This possibly works with a resulting over-estimation of the penalty, in a way depending on the unknown function $\psi$. The alternative solution is to replace the unknown quantity by an estimator (rather than a bound), and to prove that the estimator of $\psi$ built with this random penalty keeps the adaptation property of the theoretical penalized estimator. This is described in the following theorem whose proof is omitted since it is quite the same as the proof of Theorem 3.4 in Brunel and Comte (2005).

**Theorem 3.2.** *Assume that the assumptions of Theorem 3.1 are satisfied. Consider the estimator $\hat{\psi}_{\hat{m}}$ with $\hat{m}$ defined by*

$$\hat{m} = \arg \min_{m \in \mathcal{M}_n} [\gamma_n^\psi(\hat{\psi}_m) + \widehat{\mathrm{pen}}^\psi(m)]$$

and

$$\widehat{\text{pen}}^{\psi}(m) = \kappa \Phi_0^2 \left( \frac{1}{n} \sum_{i=1}^{n} \delta_i \right) \frac{D_m}{n}$$

where $\kappa$ is a universal constant, then $\hat{\psi}_{\hat{m}}$ satisfies

$$\mathbb{E}(\|\hat{\psi}_{\hat{m}} - \psi\|^2) \le \inf_{m \in \mathcal{M}_n} K_0 \left[ \|\psi - \psi_m\|^2 + \Phi_0^2 \left( \int_0^1 \psi(x) dx \right) \frac{D_m}{n} \right] + \frac{K}{n}, \quad (3.8)$$

where $K_0$ is a universal constant and $K$ depends on $\psi$, $\Phi_0$.

In particular, we can derive quite straightforwardly from results as Theorem 3.2 adaptation results to unknown smoothness:

**Proposition 3.2.** *Consider the collection of models* [T], [DP] *or* [W], *with $r > \alpha_\psi > 0$. Assume that an estimator $\tilde{\psi}$ of $\psi$ satisfies inequality (3.8) in Theorem 3.2 (respectively inequality (3.7) in Theorem 3.1). Let $L > 0$. Then*

$$\left( \sup_{\psi \in \mathbb{B}_{\alpha_\psi,2,\infty}(L)} \mathbb{E} \|\psi - \tilde{\psi}\|^2 \right)^{\frac{1}{2}} \le C(\alpha_\psi, L) n^{-\frac{\alpha_\psi}{2\alpha_\psi + 1}} \quad (3.9)$$

*where $\mathbb{B}_{\alpha_\psi,2,\infty}(L) = \{t \in \mathcal{B}_{\alpha_\psi,2,\infty}([0,1]), |t|_{\alpha_\psi,2} \le L\}$ where $C(\alpha_\psi, L)$ is a constant depending on $\alpha_\psi, L$ and also on $\psi$, $\Phi_0$.*

### 3.3. Application to the estimation of the distribution function F

Consider now the first estimator of $F$, given by (2.7).

A simple case study allows to see that if $\tilde{\psi}(x)/\tilde{g}(x) < 0$ or $\tilde{\psi}(x)/\tilde{g}(x) > 1$, then $|\tilde{F}(x) - F(x)| \le |\tilde{\psi}(x)/\tilde{g}(x) - F(x)|$, and thus the inequality $|\tilde{F}(x) - F(x)| \le |\tilde{\psi}(x)/\tilde{g}(x) - F(x)|$ holds for any $x$. Also, our definition implies that $|\tilde{F}(x) - F(x)| \le 1$, for any $x$. Moreover, to exploit [A1], we define

$$\Omega_g = \{\omega : \tilde{g}(x) - g(x) > -g_0/2, \forall x \in [0,1]\}.$$

Then, the following bounds are obtained:

$$\begin{aligned}
\|\tilde{F} - F\|^2 &= \int_0^1 (\tilde{F}(x) - F(x))^2 dx \\
&= \int_0^1 (\tilde{F}(x) - F(x))^2 dx \mathbb{I}_{\Omega_g} + \int_0^1 (\tilde{F}(x) - F(x))^2 dx \mathbb{I}_{\Omega_g^c} \\
&\le \int_0^1 \left( \frac{\tilde{\psi}(x)}{\tilde{g}(x)} - \frac{\psi(x)}{g(x)} \right)^2 dx \mathbb{I}_{\Omega_g} + \int_0^1 dx \mathbb{I}_{\Omega_g^c}.
\end{aligned}$$

Thus the first term can de decomposed as follows

$$\tilde{F} - F = \frac{\tilde{\psi} - \psi}{\tilde{g}} + F \left( \frac{g - \tilde{g}}{\tilde{g}} \right)$$

and since $\tilde{g}(x) \geq g_0/2$ on $\Omega_g$, this yields

$$\int_0^1 \left(\frac{\tilde{\psi}(x)}{\tilde{g}(x)} - \frac{\psi(x)}{g(x)}\right)^2 dx \mathbb{I}_{\Omega_g} \leq \left(\frac{2}{g_0}\right)^2 \left(\|\hat{\psi}_{\hat{m}} - \psi\|^2 + \|\tilde{g} - g\|^2\right).$$

For the second, taking the expectation, we use the following Lemma:

**Lemma 3.1.** *Assume that $g \in \mathcal{B}_{\alpha_g,2,\infty}([0,1])$ for some $\alpha_g > 1/2$ and consider a collection of spaces $S_m$ such that $\log(n) \leq D_m \leq \sqrt{n}$. Then, under Assumptions [$A_1$] and ($\mathcal{H}_2$), there exists a constant $C$ such that*

$$\mathbb{P}(\Omega_g^c) \leq \mathbb{P}\left(\|\tilde{g} - g\|_\infty > g_0/2\right) \leq \frac{C}{n}. \tag{3.10}$$

Finally, by gathering the bounds, we obtain the following proposition:

**Proposition 3.3.** *Under the assumptions of Lemma 3.1,*

$$\mathbb{E}\|\tilde{F} - F\|^2 \leq \frac{2^4}{g_0^2}\left(\mathbb{E}\|\tilde{\psi} - \psi\|^2 + \mathbb{E}\|\tilde{g} - g\|^2\right) + \frac{C(g_0, \|\psi\|)}{n}, \tag{3.11}$$

*where $C(g_0, \|\psi\|)$ is a constant depending on $g_0$ and $\|\psi\|$.*

From Inequality (3.11), we easily deduce by using results (3.7) or (3.8) that $\tilde{F}$ is an adaptive estimator of $F$ if the functions $g$ and $\psi$ have the same regularity $\alpha = \alpha_g = \alpha_\psi$. Here we can state the following result:

**Proposition 3.4.** *Assume that $g \in \mathcal{B}_{\alpha_g,2,\infty}([0,1])$ and that $\psi \in \mathcal{B}_{\alpha_\psi,2,\infty}([0,1])$. Consider the collection of models [T], [DP] or [W], with dimensions $\log(n) \leq D_m \leq \sqrt{n}$ and with $r > \alpha_F = \alpha_\psi = \alpha_g > 1/2$. Let $\tilde{F}$ the estimator defined by (2.7) and let $L > 0$. Then*

$$\left(\sup_{F \in \mathbb{B}_{\alpha_F,2,\infty}(L)} \mathbb{E}\|F - \tilde{F}\|^2\right)^{\frac{1}{2}} \leq C(\alpha_F, L) n^{-\frac{\alpha_F}{2\alpha_F + 1}} \tag{3.12}$$

*where $\mathbb{B}_{\alpha_F,2,\infty}(L) = \{t \in \mathcal{B}_{\alpha_F,2,\infty}([0,1]), |t|_{\alpha_F,2} \leq L\}$ where $C(\alpha_F, L)$ is a constant depending on $\alpha_F, L$ and also on $\psi$, $\Phi_0$ and $g_0$.*

Note that Theorem 2 in Yang (2000) shows that the rate in the sup-norm over a compact is of order $O((\ln n/n)^{(1+\alpha_f)/(3+2\alpha_f)})$ a.s. where $\alpha_f$ stands for the regularity of the density function $f$ (that is $\alpha_F = \alpha_f + 1$).

If the index of regularity of $F$, $\alpha_F$, is greater than the index of regularity of $\psi = Fg$, $\alpha_\psi$, then the asymptotic rate of the estimator $\tilde{F}$ is given by $n^{-\alpha_\psi/(1+2\alpha_\psi)}$ instead of the optimal one $n^{-\alpha_F/(1+2\alpha_F)}$. This is the reason why we propose another contrast to estimate directly $F$.

## 4. Study of the mean square estimator

In this section, we study the mean square estimator of $F$ from (2.9) and its adaptive version. In this context, we define the empirical norm $\|\cdot\|_n$ as follows: for $t \in S_m$, $\|t\|_n^2 = (1/n) \sum_{i=1}^n t^2(U_i)$. It is a natural norm in regression problems, and under [A1], it is equivalent in mean to the standard Lebesgue integrated $\mathbb{L}^2$-norm, i.e., under [A1]:

$$\forall t \in S_m, \ g_0 \|t\|^2 \leq \mathbb{E}(\|t\|_n^2) = \int t^2(x) g(x) dx \leq g_1 \|t\|^2.$$

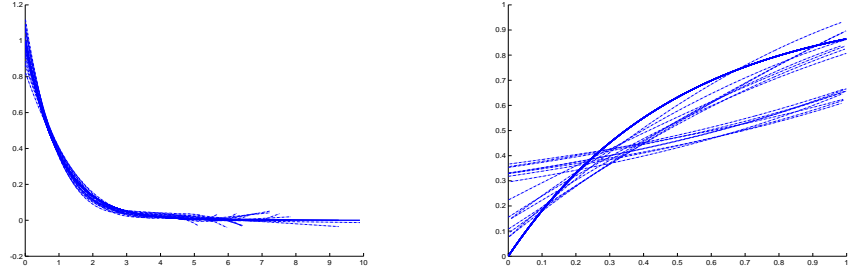

FIG 1. *Plot of 15 Quotient estimators (left: density estimator $\tilde{g}$ of $g$, right: Quotient estimator $\tilde{F}$ of $F$) for Model 4 with $n = 500$.*

Then, the mean-square contrast defined by (2.8) can be decomposed as follows:

$$\gamma_n^{MS}(t) - \gamma_n^{MS}(s) = \|t - F\|_n^2 - \|s - F\|_n^2 - 2\nu_n^{MS}(t - s) \tag{4.1}$$

where $\nu_n^{MS}(.)$ is defined by:

$$\nu_n^{MS}(t) = \frac{1}{n} \sum_{i=1}^n (\delta_i - F(U_i)) t(U_i) \tag{4.2}$$

which is a centered process since $\mathbb{E}(\delta | U = u) = F(u)$.

In this case, we obtain the following result for the penalized estimator:

**Theorem 4.1.** *Consider the collections of models [T] with $N_n \leq \sqrt{n}/\ln(n)$ or [DP] or [W] with $N_n \leq n/\ln^2(n)$. Let $\hat{F}_{\hat{m}_0}$ be defined by (2.10), with*

$$\text{pen}^{MS}(m) \geq \kappa_0 \frac{D_m}{n}.$$

*Then,*

$$\mathbb{E}(\|\hat{F}_{\hat{m}_0} - F\|_n^2) \leq C \inf_{m \in \mathcal{M}_n} (\|F_m - F\|_n^2 + \text{pen}^{MS}(m)) + C' \frac{1}{n} \tag{4.3}$$

*where $F_m$ stands for the orthogonal projection of $F$ on $S_m$ and $C$ and $C'$ are constants depending on $\Phi_0$ and $g$.*

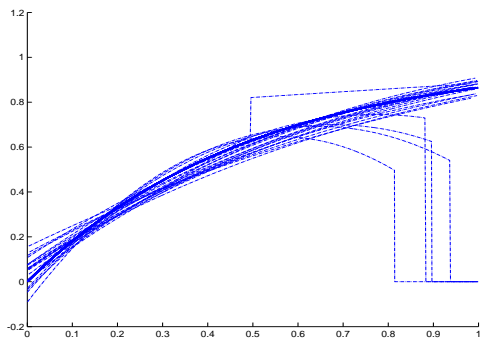

FIG 2. *Plot of 15 Regression estimators $\hat{F}_{\hat{m}_0}$ for Model 4 with $n = 500$.*

Note that the computation of the estimator may be more tedious in practice than the quotient one, but the result is obtained directly for the estimator of $F$, without any regularity condition on $\psi$. As a consequence, we obtain here a rate only depending on the regularity of $F$, and we can state the following result:

**Proposition 4.1.** *Consider the collection of models* [T] *with $\alpha_F > 1/2$ or* [DP] *or* [W]*, with $r > \alpha_F > 0$ and the estimator $\hat{F}_{\hat{m}_0}$ defined by (2.10)–(2.11). Let $L > 0$. Then*

$$\left( \sup_{F \in \mathbb{B}_{\alpha_F,2,\infty}(L)} \mathbb{E}\|F - \hat{F}_{\hat{m}_0}\|_n^2 \right)^{\frac{1}{2}} \leq C(\alpha_F, L) n^{-\frac{\alpha_F}{2\alpha_F + 1}} \tag{4.4}$$

*where $\mathbb{B}_{\alpha_F,2,\infty}(L) = \{t \in \mathcal{B}_{\alpha_F,2,\infty}([0,1]), |t|_{\alpha_F,2} \leq L\}$ where $C(\alpha_F, L)$ is a constant depending on $\alpha_F, L$ and also on $F$, $\Phi_0$ and $g_0$.*

## 5. Simulations

We consider the regular collection [DP] (see Section 2.3) with degrees less than $r_{max}$ on a subdivision $[j/2^p, (j+1)/2^p[$. The density and regression algorithms minimize the contrast and select the approximation space in the sense that the integers $p$ and $r$ are selected such that $2^p(r+1) \leq N_n \leq n/\log^2(n)$ and $r \in \{0, 1, \ldots, r_{max}\}$. Note that the degree $r$ is global in the sense that it is the same on all the intervals of the subdivision. We take $r_{max} = 9$ in practice. Moreover, additive (but negligible) correcting terms are classically involved in the penalty (see Comte and Rozenholc (2004)). Such terms avoid under-penalization and are in accordance with the fact that the theorems provide lower bounds for the penalty. As the correcting terms are asymptotically negligible, they do not affect the rate of convergence. The constants in the penalty are taken equal to 4. Finally, for $m = (p, r)$, the penalties are proportional to $2^p(r + 1 + \log^{2.5}(r+1))$ with proportionality factor $\kappa = 4$ for the estimation of $g$ and $F$ and a factor $(4/n) \sum_{i=1}^{n} \delta_i$ for the estimation of $\psi$. Most programs are available on Yves Rozenholc's web page http://www.math-info.univ-paris5.fr/~rozen/.

TABLE 1
*Monte-Carlo results for the MSE ($\times 10^{-2}$) of the quotient, regression and the NPMLE estimators of the cdf $F$, for $J = 200$ or $500$ sample replications.*

|  | Quotient est. | | | | Regression est. | | | |
|---|---|---|---|---|---|---|---|---|
| $n$ | 60 | 200 | 500 | 1000 | 60 | 200 | 500 | 1000 |
| model 1 | 1.75 | 0.49 | 0.16 | 0.08 | 0.50 | 0.13 | 0.05 | 0.03 |
| model 2 | 2.73 | 0.83 | 0.38 | 0.22 | 10.3 | 0.89 | 0.23 | 0.01 |
| model 3 | 1.99 | 0.64 | 0.28 | 0.09 | 1.13 | 0.39 | 0.07 | 0.03 |
| model 4 | 5.96 | 3.99 | 1.86 | 0.36 | 7.40 | 2.40 | 0.48 | 0.19 |
| model 5 | 1.89 | 0.79 | 0.37 | 0.18 | 0.76 | 0.27 | 0.11 | 0.07 |
|  | Birgé's NPMLE | | | | Groeneboom's NPMLE | | | |
| $n$ | 60 | 200 | 500 | 1000 | 60 | 200 | 500 | 1000 |
| model 1 | 1.72 | 0.75 | 0.41 | 0.25 | 1.94 | 0.75 | 0.37 | 0.24 |
| model 2 | 2.15 | 0.75 | 0.51 | 0.27 | 1.97 | 0.76 | 0.40 | 0.23 |
| model 3 | 1.52 | 0.76 | 0.39 | 0.25 | 2.15 | 0.80 | 0.40 | 0.22 |
| model 4 | 2.90 | 1.11 | 0.66 | 0.38 | 2.04 | 0.82 | 0.43 | 0.26 |
| model 5 | 1.20 | 0.93 | 0.32 | 0.27 | 0.93 | 0.38 | 0.20 | 0.12 |

Now, let us describe the simulated models. Remember that the distribution of $\delta$ given $U = u$ is a Bernoulli variable with parameter $F(u)$. We consider the following models for generating data:

Model 1. *Uniform distribution $F$*: $U \sim \mathcal{U}(0, 1)$ and $\delta \sim \mathcal{B}(1, U)$
Model 2. *$\chi^2$-distribution $F$*: $U \sim \mathcal{U}(0, 1)$ and $\delta \sim \mathcal{B}(1, F_{\chi_1^2}(U))$
Model 3. *Quadratic distribution $F$*: $U \sim \mathcal{U}(0, 1)$ and $\delta \sim \mathcal{B}(1, U^2)$
Model 4. *Exponential distribution $F$*: $U \sim \gamma(1, \lambda)$ and $\delta \sim \mathcal{B}(1, 1 - e^{-\mu U})$ with $\lambda = 1, \mu = 0.5$.
Model 5. *Beta distribution (S-shape) $F$*: $U \sim \beta(4, 6)$ and $\delta \sim \mathcal{B}(1, F_{\beta(4,8)}(U))$ where $F_{\beta(\alpha,\beta)}$ is the cdf of a Beta distribution of parameter $(\alpha, \beta)$.

To study the quality of each estimation procedure and to compare them, we compute over $J$ sample replications of size $n = 60, 200, 500$ and $1000$ the mean squared errors (MSE) over the sample points $u_1, \ldots, u_K$ falling in $[a, b]$:

$$\text{MSE}_j = \frac{(b-a)}{K} \sum_{k=1}^{K} [F(u_k) - \hat{F}_j(u_k)]^2$$

where $\hat{F}_j$ stands for the (adaptive) quotient estimator $\tilde{F}$ or for the penalized regression estimator $\hat{F}_{\hat{m}_0}$ or the benchmark NPMLEs, computed over the $j$th sample replication for $j = 1, \ldots, J$. For the small sample sizes $n = 60$ and $n = 200$, the average values are obtained with $J = 500$ repetitions while for large samples ($n = 500$ and $n = 1000$), only $J = 200$ replications are performed. To avoid boundary effects due to the sparsity of the observations at the end of the interval particularly for models 2, 4, and 5 the $\text{MSE}_j$'s are truncated for each replication in the sense that we include in the mean only the $u_k$ less than a given quantile value: $\mathbb{P}(X \leq 1) = 0.68$ for model 2, $\mathbb{P}(X \leq 1) = 0.86$ for model 4 and $\mathbb{P}(X \leq 0.5) = 0.89$ for model 5; thus, the $\text{MSE}_j$ are computed over $[a, b]$ with $a = 0$ and $b = 1$ from model 1 to 4, and $b = 0.5$ for model 5. Therefore, the MSE's given in Table 1 stand effectively for the truncated arithmetic means of the $\text{MSE}_j$'s.

As we can see from results in Table 1, the regression estimator is always better than the quotient and the NPMLE estimators for large samples. However, for small sample sizes, the quotient estimator can behave as well as and even better than the regression one, see models 2 and 4 for $n = 60, 200$. The same remark holds for both Birgé's and Groeneboom's NPMLE. These last estimators have the advantage to be very easy to compute. As a counterpart they look like step functions whatever the regularity of the function is. Moreover, Birgé's estimator is not adaptive since we have to choose the number of partition cells ($D = 5$ cells for a sample size $n = 60, 200$ and $D = 10$ cells for $n = 500, 1000$), see Figure 3. Note also that, the density estimator $\tilde{g}$ of $g$ is a very attractive estimator by itself as shown in Figure 1. In some cases and particularly for model 4, see Figure 1, the quotient mechanism works wrong even if the density estimator is very performing. Figure 1 (right) shows that near than half of the curves do not give the good shape. This is a drawback of quotient strategies which do not have good robustness properties. Regression estimators (see Figure 2) are much more stable. In Figure 3, we give an illustration of all the compared estimators for small ($n = 60$) and large ($n = 1000$) samples and we can see that our adaptive estimators behave as well as and often better than the benchmark NPMLEs.

**Concluding remarks.** Globally it appears that the regression estimation is better than the quotient estimator, from both theoretical and empirical points of view. The latter can be better than the former only for small sample experiments. The two density estimators involved in the quotient are nevertheless easy to compute, and empirically very good. It is thus interesting to see that the estimation algorithms give very good results. Nevertheless, even for well estimated numerator and denominator, the ratio is less satisfactory than the direct regression estimator.

## 6. Proofs

### *6.1. Talagrand's Inequality*

The following version of Talagrand's Inequality (see Talagrand (1996)) is very useful in most of the proofs:

**Lemma 6.1.** *Let $Z_1, \ldots, Z_n$ be i.i.d. random variables and $\nu_n(\ell)$ be defined by $\nu_n(\ell) = (1/n) \sum_{i=1}^{n} [\ell(Z_i) - \mathbb{E}(\ell(Z_i))]$ for $\ell$ belonging to a countable class $\mathcal{L}$ of uniformly bounded measurable functions. Then there exists a universal constant $K_0$ such that, for any positive $\eta$, $\lambda$*

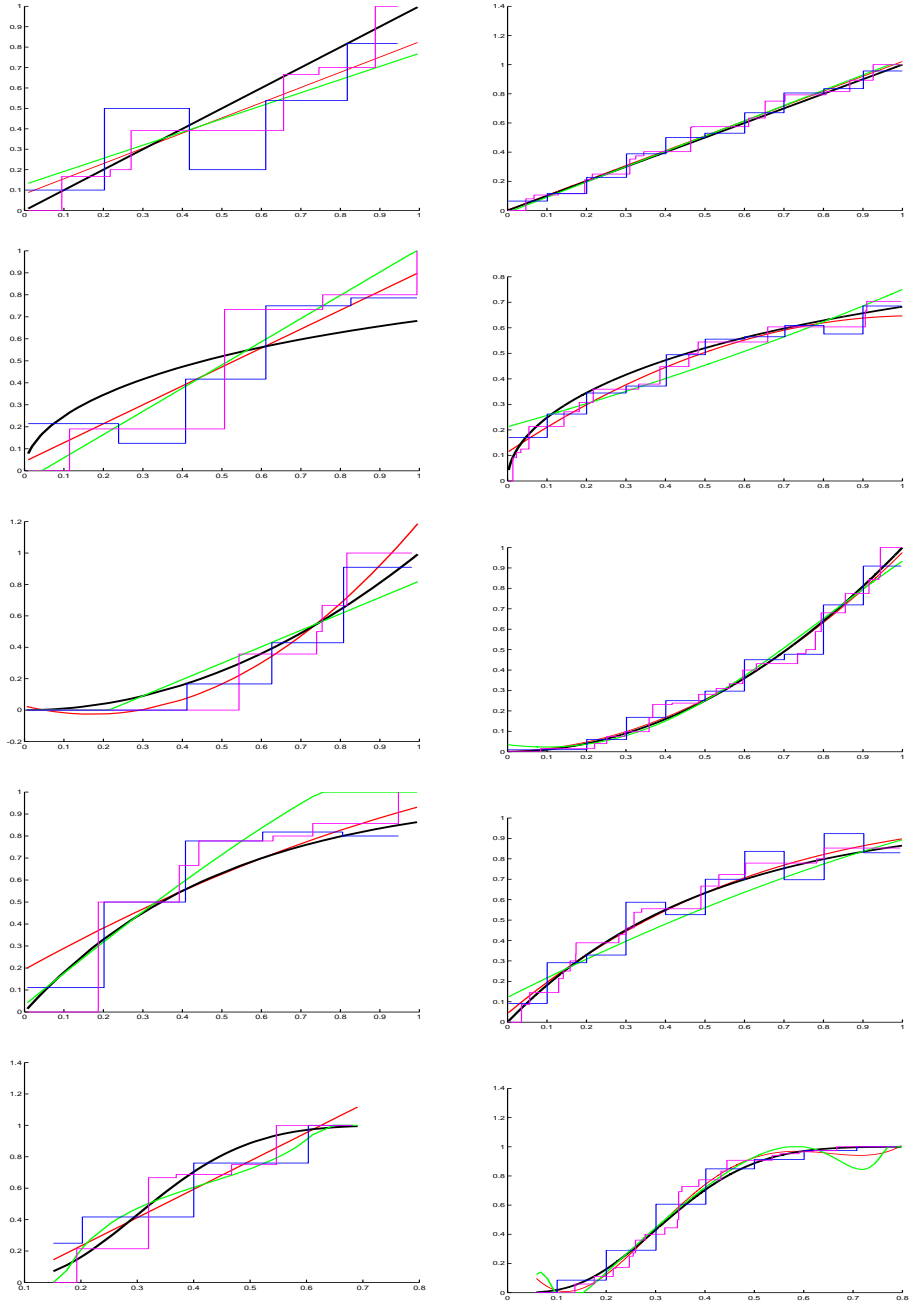

FIG 3. *True cdf (black thick line), Regression estimator (red), Quotient estimator (green), Birgé's NPMLE (blue) and Groeneboom's NPMLE (magenta) for Model 1 to 5 (top to bottom) with $n = 60$ (left) and $n = 1000$ (right).*

$$\mathbb{P}(\sup_{\ell \in \mathcal{L}} |\nu_n(\ell)| \geq (1+\eta)H + \lambda) \leq 3\exp\left[-K_0 n\left(\frac{\lambda^2}{v} \wedge \frac{(\eta \wedge 1)\lambda}{M_1}\right)\right]. \quad (6.1)$$

Moreover, there exists a universal constant $K_1$ such that for $\epsilon > 0$

$$\mathbb{E}\left[\sup_{\ell \in \mathcal{L}} |\nu_n(\ell)|^2 - 2(1+2\epsilon)H^2\right]_+ \leq \frac{6}{K_1}\left(\frac{v}{n}e^{-K_1\epsilon\frac{nH^2}{v}} + \frac{8M_1^2}{K_1 n^2 C^2(\epsilon)}e^{-\frac{K_1 C(\epsilon)\sqrt{\epsilon}}{\sqrt{2}}\frac{nH}{M_1}}\right), \quad (6.2)$$

with $C(\epsilon) = (\sqrt{1+\epsilon} - 1) \wedge 1$, and where

$$\sup_{\ell \in \mathcal{L}} \|\ell\|_\infty \leq M_1, \quad \mathbb{E}\left(\sup_{\ell \in \mathcal{L}} |\nu_n(\ell)|\right) \leq H, \quad \sup_{\ell \in \mathcal{L}} \text{Var}(\ell(X_1)) \leq v.$$

Note that (6.1) is also given in Birgé and Massart (1998), Corollary 2. In both cases, usual density arguments allow to take instead of the class $\mathcal{L}$ a unit ball in a finite dimension space of functions.

### 6.2. Proof of Lemma 3.1

Let us write
$$\|\tilde{g} - g\|_\infty \leq \|g - g_{\hat{m}_g}\|_\infty + \|g_{\hat{m}_g} - \hat{g}_{\hat{m}_g}\|_\infty$$

with $\tilde{g} = \hat{g}_{\hat{m}_g}$ defined by (2.4). As $g$ belongs to some Besov space $\mathcal{B}_{\alpha_g,2,\infty}([0,1])$ with $\alpha_g > 1/2$ and as $\mathcal{B}_{\alpha_g,2,\infty}([0,1]) \subset \mathcal{B}_{\alpha_g-1/2,\infty,\infty}([0,1])$ then, Lemma 12 in Barron et al. (1999) gives (with the restriction $D_m \geq \log(n), \forall m$):

$$\|g - g_{\hat{m}_g}\|_\infty \leq CD_{\hat{m}_g}^{-(\alpha_g-1/2)} \leq C(\log n)^{-(\alpha_g-1/2)}.$$

Thus, $\|g - g_{\hat{m}_g}\|_\infty$ decreases to 0 as $n$ goes to $\infty$ and for some integer $n_0$ large enough, we have for $n \geq n_0$,

$$\mathbb{P}(\|\tilde{g} - g\|_\infty > g_0/2) \leq \mathbb{P}(\|g_{\hat{m}_g} - \hat{g}_{\hat{m}_g}\|_\infty > g_0/4)$$

Now, $\|g_{\hat{m}_g} - \hat{g}_{\hat{m}_g}\|_\infty \leq \Phi_0\sqrt{D_{\hat{m}_g}}\|g_{\hat{m}_g} - \hat{g}_{\hat{m}_g}\|$ and $\|g_{\hat{m}_g} - \hat{g}_{\hat{m}_g}\|^2 = \sum_{\lambda \in \Lambda_{\hat{m}_g}} \nu_{n,g}^2(\varphi_\lambda) = \sup_{t \in B_{\hat{m}_g}} |\nu_{n,g}^2(t)|$. This implies

$$\begin{aligned}\mathbb{P}(\|\tilde{g} - g\|_\infty > g_0/2) &\leq \mathbb{P}\left(\sup_{t \in B_{\hat{m}_g}} |\nu_{n,g}(t)| > \frac{g_0}{4\Phi_0\sqrt{D_{\hat{m}_g}}}\right) \\ &\leq \sum_{m \in \mathcal{M}_n} \mathbb{P}\left(\sup_{t \in B_m} |\nu_{n,g}(t)| > \frac{g_0}{4\Phi_0\sqrt{D_m}}\right) \quad (6.3)\end{aligned}$$

We apply Inequality (6.1) to the class of functions $\mathcal{L} = B_m(0,1)$ by taking $\ell = t - \mathbb{E}(t(U_1))$, with

$$\sup_{t \in B_m(0,1)} \|t\|_\infty \leq \Phi_0\sqrt{D_m} := M_1,$$
$$\sup_{t \in B_m(0,1)} \text{Var}(t(U_1)) \leq \sup_{t \in B_m(0,1)} \int_0^1 t^2(u)g(u)du \leq g_1 := v$$

and $\mathbb{E}(\sup_{t\in B_m(0,1)} \nu_{n,g}^2(t)) = (1/n)\sum_{\lambda\in\Lambda_m} \text{Var}(\varphi_\lambda(U_1)) \leq \Phi_0^2 D_m/n := H^2$. By choosing $\eta = 1$ and $\lambda = g_0/(8\Phi_0\sqrt{D_m})$ and if $2E + \lambda \leq g_0/(4\Phi_0\sqrt{D_m})$, we obtain from (6.3):

$$\mathbb{P}(\|\tilde{g} - g\|_\infty > g_0/2) \leq \sum_{m\in\mathcal{M}_n} 3\exp\left[-K_1 n\left(\frac{\lambda^2}{g_1} \wedge \frac{\lambda}{\Phi_0\sqrt{D_m}}\right)\right]$$

$$\leq \sum_{m\in\mathcal{M}_n} 3\exp\left[-K_1 C_1 n/D_m\right]$$

with $C_1 = \left(\frac{g_0^2}{64g_1\Phi_0^2} \wedge \frac{g_0}{8\Phi_0^2}\right)$, if we ensure that $2E+\lambda \leq g_0/(4\Phi_0\sqrt{D_m})$. But with $E = \Phi_0\sqrt{D_m/n}$, this is verified if $D_m \leq [g_0/(16\Phi_0^2)]\sqrt{n}$. Thus, we can deduce that

$$\mathbb{P}(\|\tilde{g} - g\|_\infty > g_0/2) \leq 3|\mathcal{M}_n|\exp\left[-K_1 C_1'\sqrt{n}\right]$$

with $C_1' = \left(\frac{g_0^2}{64g_1\Phi_0^2} \wedge \frac{g_0}{8\Phi_0^2}\right)/[g_0/(16\Phi_0^2)]$. Finally, since $|\mathcal{M}_n| \leq n$, if $D_m \leq (K_1 C_1)n/(2\ln(n))$ then $\mathbb{P}(\|\tilde{g} - g\|_\infty > g_0/2) \leq 3/n$ and this concludes the proof. $\square$

### 6.3. Proof of Theorem 3.1

#### 6.3.1. Proof of a preliminary Lemma

First, we prove the following lemma:

**Lemma 6.2.** *Assume that $(\mathcal{H}_1)$ and $(\mathcal{H}_2)$ are fulfilled and denote by $B_{m,m'}(0,1) = \{t \in S_m + S_{m'}, \ \|t\| = 1\}$. Let $\nu_n(\ell_t)$ be defined by (3.5) and*

$$\ell_t(u,\delta) = \delta t(u), \tag{6.4}$$

*then for $\epsilon > 0$*

$$\mathbb{E}\left(\sup_{t\in B_{m,m'}(0,1)} \nu_n^2(\ell_t) - p^\psi(m,m')\right)_+ \leq \frac{\kappa_1}{n}\left(e^{-\kappa_2\epsilon(D_m+D_{m'})} + \frac{e^{-\kappa_3\epsilon^{3/2}\sqrt{n}}}{C(\epsilon)^2}\right), \tag{6.5}$$

*with $p^\psi(m,m') = 2(1+2\epsilon)\Phi_0^2 \int_0^1 \psi(x)dx\,(D_m + D_{m'})/n$ and $C(\epsilon) = (\sqrt{1+\epsilon} - 1) \wedge 1$. The constants $\kappa_i$ for $i = 1,2,3$ depend on $\Phi_0, \psi$ and $F$.*

We apply Talagrand's inequality (6.2) by taking $Z_i = (U_i, \delta_i)$ for $i = 1,\ldots,n$ and $\ell(u,\delta) = \ell_t(u,\delta)$. Usual density arguments show that this result can be applied to the class of functions $\mathcal{L} = \{\ell_t, t \in B_{m,m'}(0,1)\}$. Then we find for the present empirical process the following bounds:

$$\sup_{\ell\in\mathcal{L}} \|\ell\|_\infty = \sup_{t\in B_{m,m'}(0,1)} \|\ell_t\|_\infty \leq \Phi_0\sqrt{D(m,m')} := M_1$$

with $D(m, m')$ denoting the dimension of $S_m + S_{m'}$. Then

$$\sup_{\ell \in \mathcal{L}} \text{Var}(\ell(U_1, \delta_1)) = \sup_{t \in B_{m,m'}(0,1)} \text{Var}(\ell_t(U_1, \delta_1)) = \sup_{t \in B_{m,m'}(0,1)} \mathbb{E}(\delta_1 t^2(U_1))$$

$$= \sup_{t \in B_{m,m'}(0,1)} \int_0^1 t^2(u)\psi(u)du \leq g_1 := v.$$

Lastly,

$$\mathbb{E}\left(\sup_{\ell \in \mathcal{L}} \nu_n^2(\ell)\right) = \mathbb{E}\left(\sup_{t \in B_{m,m'}(0,1)} \nu_n^2(\ell_t)\right) \leq \sum_{\lambda \in \Lambda_{m,m'}} \frac{1}{n}\text{Var}(\delta_1 \varphi_\lambda(U_1))$$

$$\leq \frac{\Phi_0^2 D(m, m')}{n} \int_0^1 \psi(x)dx = C_1 \frac{D(m, m')}{n} := H^2.$$

with the natural notation $\Lambda_{m,m'} = \Lambda_m \cup \Lambda_{m'}$. Then it follows from (6.2) that

$$\mathbb{E}\left(\sup_{t \in B_{m,m'}(0,1)} \nu_n^2(\ell_t) - p^\psi(m, m')\right) \leq \kappa_1 \left(\frac{1}{n} e^{-\kappa_2 \epsilon D(m')} + \frac{1}{nC^2(\epsilon)} e^{-\kappa_3 \epsilon^{3/2} \sqrt{n}}\right),$$

where $\kappa_i$ for $i = 1, 2, 3$ are constant depending on $K_1$ and $C_1$ and $p^\psi(m, m') = 2(1 + 2\epsilon)C_1(D_m + D_{m'})/n$. □

### 6.3.2. Proof of Theorem 3.1

It follows from the definition of $\hat{\psi}_{\hat{m}}$ that: $\forall m \in \mathcal{M}_n$,

$$\gamma_n^\psi(\hat{\psi}_{\hat{m}}) + \text{pen}^\psi(\hat{m}) \leq \gamma_n^\psi(\psi_m) + \text{pen}^\psi(m). \tag{6.6}$$

Then by using decomposition (3.6), it follows from (6.6) and from the definition of the process $\nu_n(\ell_t)$ given by (3.5) that:

$$\|\hat{\psi}_{\hat{m}} - \psi\|^2 \leq \|\psi_m - \psi\|^2 + 2\nu_n(\ell_{\hat{\psi}_{\hat{m}} - \psi}) + \text{pen}^\psi(m) - \text{pen}^\psi(\hat{m})$$

$$\leq \|\psi_m - \psi\|^2 + \frac{1}{4}\|\hat{\psi}_{\hat{m}} - \psi_m\|^2 + 4\sup_{t \in B_{m,\hat{m}}(0,1)} \nu_n^2(\ell_t)$$

$$+ \text{pen}^\psi(m) - \text{pen}^\psi(\hat{m})$$

where we recall that $B_{m,\hat{m}}(0,1) = \{t \in S_m + S_{\hat{m}} \ / \ \|t\| \leq 1\}$. Note that the norm connection as described by (2.3) still holds for any element $t$ of $S_m + S_{m'}$ as follows: $\|t\|_\infty \leq \Phi_0 \max(D_m, D_{m'})\|t\|$. Indeed, under $(\mathcal{H}_2)$, we restrict our attention to nested collections of models, so that $S_m + S_{\hat{m}}$ is equal to the larger of the two spaces. For a fixed integer $m$, we denote by $D(m')$ the dimension of $S_m + S_{m'}$, for all $m' \in \mathcal{M}_n$. Note that $D(m') = \max(D_m, D_{m'}) \leq D_m + D_{m'}$.

Let $p^\psi(m, m')$ be defined as in Lemma 6.2. Then $\forall m \in \mathcal{M}_n$,

$$\frac{1}{2}\|\hat{\psi}_{\hat{m}} - \psi\|^2 \leq \frac{3}{2}\|\psi - \psi_m\|^2 + 2\,\text{pen}^\psi(m) + 8\left(\sup_{t \in B_{m,\hat{m}}(0,1)} \nu_n^2(\ell_t) - p^\psi(m, \hat{m})\right)$$

$$+ 8p^\psi(m, \hat{m}) + \text{pen}(m) - \text{pen}(\hat{m}).$$

Now, note first that

$$\mathbb{E}\Big(\sup_{t\in B_{m,\hat{m}}(0,1)} \nu_n^2(\ell_t) - p^\psi(m,\hat{m})\Big)_+ \leq \sum_{m'\in\mathcal{M}_n} \mathbb{E}\Big(\sup_{t\in B_{m,m'}(0,1)} \nu_n^2(\ell_t) - p^\psi(m,m')\Big)_+. \tag{6.7}$$

Moreover it follows from Lemma 6.2 that

$$\sum_{m'\in\mathcal{M}_n} \mathbb{E}\Big(\sup_{t\in B_{m,m'}(0,1)} \nu_n^2(\ell_t) - p^\psi(m,m')\Big)_+ \leq \kappa_1 \left(\frac{\Sigma(m')}{n} + \frac{|\mathcal{M}_n|}{n}e^{-\kappa_3\epsilon^{3/2}\sqrt{n}}\right)$$

where $\Sigma(m') = \sum_{m'\in\mathcal{M}_n} e^{-\kappa_2\epsilon D(m')}$. Then by taking $\epsilon = 1/2$ and assuming that $|\mathcal{M}_n| \leq n$ and since, under $(\mathcal{H}_2)$, $\sum_{m\in\mathcal{M}_n} e^{-aD_m} \leq \sum_{k=1}^n e^{-ka} \leq \Sigma(a) < +\infty, \forall a > 0$, this leads to the bound

$$\mathbb{E}\Big(\sup_{t\in B_{m,\hat{m}}(0,1)} \nu_n^2(\ell_t) - p^\psi(m,m')\Big)_+ \leq \frac{C}{n}.$$

Therefore, we have the following result, which proves the theorem: $\forall m \in \mathcal{M}_n$,

$$\mathbb{E}(\|\hat{\psi}_{\hat{m}} - \psi\|^2) \leq 3\|\psi - \psi_m\|^2 + 4\text{pen}^\psi(m) + \frac{C}{n} + 2\mathbb{E}\left(8p^\psi(m,\hat{m}) - \text{pen}(m) - \text{pen}(\hat{m})\right).$$

Therefore by using the definition of $p^\psi(m,m')$ in Lemma 6.2, we choose

$$\text{pen}^\psi(m) \geq 16(1+2\epsilon) \int_0^1 \psi(x)dx \frac{D_m}{n}.$$

This ensures that $\forall m, m', 8p^\psi(m,m') \leq \text{pen}(m) + \text{pen}(m')$ and yields to (6.7). $\square$

### 6.4. Proof of Theorem 3.2

We start by writing that, $\forall m \in \mathcal{M}_n$,

$$\gamma_n(\hat{F}_{\hat{m}_0}) + \text{pen}^{\text{MS}}(\hat{m}_0) \leq \gamma_n(F_m) + \text{pen}^{\text{MS}}(m)$$

and by using the decomposition (4.1). It follows that

$$\|\hat{F}_{\hat{m}_0} - F\|_n^2 \leq \|F_m - F\|_n^2 + 2\nu_n^{MS}(\hat{F}_{\hat{m}_0} - F_m) + \text{pen}^{\text{MS}}(m) - \text{pen}^{\text{MS}}(\hat{m}_0).$$

In the same way as Baraud *et al.* (2001), we introduce for $\|t\|_g^2 = \int_0^1 t^2(u)g(u)\,du$, the ball $B_{m,m'}^g(0,1) = \{t \in S_m + S_{m'}, \|t\|_g = 1\}$ and the set

$$\Omega_n = \left\{\omega, \left|\frac{\|t\|_n^2}{\|t\|_g^2} - 1\right| \leq \frac{1}{2}, \forall t \in \bigcup_{m,m'\in\mathcal{M}_n}(S_m + S_{m'}) \setminus \{0\}\right\}.$$

On the complement of $\Omega_n$, a separate study leads to the following lemma:

**Lemma 6.3.** *If $N_n \leq \sqrt{n}/\ln(n)$ for [T] or $N_n \leq n/\ln^2(n)$ for [P] or [W], then $\mathbb{P}(\Omega_n^c) \leq c/n$ and, $\mathbb{E}(\|\hat{F}_{\hat{m}_0} - F\|_n^2 \mathbb{I}_{\Omega_n^c}) \leq c'/n$, where $c$ and $c'$ are positive constants.*

*Proof of Lemma 6.3.* That $\mathbb{P}(\Omega_n^c) \leq c/n^2$ is in fact a pure property of the basis and is proved under our assumptions in Baraud. (2002). Moreover, $\|\hat{F}_{\hat{m}_0} - F\|_n^2 \leq 2(\|\hat{F}_{\hat{m}_0}\|_n^2 + \|F\|_n^2)$. Now $\|F\|_n^2 \leq 1$ and $\|\hat{F}_{\hat{m}_0}\|_n^2 = (1/n)\|\Pi_{\hat{m}_0}\delta\|_{\mathbb{R}^n}^2$ where $\delta = (\delta_1, \ldots, \delta_n)$, $\Pi_{\hat{m}}$ is the orthogonal projection in $\mathbb{R}^n$ on $\{t(U_1), \ldots, t(U_n)), t \in S_m\}$ and $\|\cdot\|_{\mathbb{R}^n}$ is the Euclidean norm in $\mathbb{R}^n$. It follows that $\|\hat{F}_{\hat{m}_0}\|_n^2 \leq (1/n)\|\delta\|_{\mathbb{R}^n}^2 = (1/n)\sum_{i=1}^n \delta_i^2 \leq 1$. Therefore $\mathbb{E}(\|\hat{F}_{\hat{m}_0} - F\|_n^2 \mathbb{I}_{\Omega_n^c}) \leq 2\mathbb{P}(\Omega_n^c) \leq c'/n$. □

Therefore, we focus on the study of the bounds on $\Omega_n$, where the inequality $\|t\|_g^2 \leq 2\|t\|_n^2$ is fulfilled. We obtain

$$\begin{aligned}\|\hat{F}_{\hat{m}_0} - F\|_n^2 \mathbb{I}_{\Omega_n} &\leq \|F_m - F\|_n^2 + \frac{1}{8}\|\hat{F}_{\hat{m}_0} - F_m\|_f^2 \mathbb{I}_{\Omega_n} + 16 \sup_{t \in B_{\hat{m}_0,m}^g(0,1)}[\nu_n^{MS}]^2(t) \\ &\quad + \text{pen}^{MS}(m) - \text{pen}^{MS}(\hat{m}_0) \\ &\leq \left(1 + \frac{1}{2}\right)\|F_m - F\|_n^2 + \frac{1}{2}\|\hat{F}_{\hat{m}_0} - F\|_n^2 \mathbb{I}_{\Omega_n} \\ &\quad + 16\left(\sup_{t \in B_{\hat{m}_0,m}^g(0,1)}[\nu_n^{MS}]^2(t) - \tilde{p}(m, \hat{m}_0)\right)_+ \\ &\quad + \text{pen}^{MS}(m) + 16\tilde{p}(m, \hat{m}_0) - \text{pen}^{MS}(\hat{m}_0).\end{aligned}$$

Let $(\bar{\varphi}_\lambda)_{\lambda \in \Lambda_{m,m'}}$ be an orthonormal basis of $S_m + S_{m'}$ for the scalar product $\langle \cdot, \cdot \rangle_g$ (built by Gramm-Schmidt orthonormalization). It is easy to see that:

$$\begin{aligned}\mathbb{E}\left(\sup_{t \in B_{m',m}^g(0,1)}[\nu_n^{MS}]^2(t)\right) &\leq \sum_{\lambda \in \Lambda_{m,m'}} \frac{1}{n}\text{Var}\left([\delta_1 - F(U_1)]\bar{\varphi}_\lambda(U_1)\right) \\ &\leq \sum_{\lambda \in \Lambda_{m,m'}} \frac{1}{n}\mathbb{E}_X\left(\int_0^1 [\mathbb{I}_{X \leq u} - F(u)]^2 \bar{\varphi}_\lambda(u)^2 g(u) du\right) \\ &\leq \frac{1}{n}\sum_{\lambda \in \Lambda_{m,m'}}\left(\int_0^1 \mathbb{E}_X[\mathbb{I}_{X \leq u} - F(u)]^2 \bar{\varphi}_\lambda^2(u) g(u) du\right) \\ &\leq \frac{1}{n}\sum_{\lambda \in \Lambda_{m,m'}}\left(\int_0^1 F(u)(1 - F(u))\bar{\varphi}_\lambda^2(u) g(u) du\right) \\ &\leq \frac{D_m \vee D_{m'}}{n}\end{aligned}$$

as $F(u)(1 - F(u)) \leq 1$. Therefore, we obtain by applying Talagrand's Inequality

$$\sum_{m' \in \mathcal{M}_n} \mathbb{E}\left(\sup_{t \in B_{m',m}^g((0,1)}[\nu_n^{MS}]^2(t) - \tilde{p}(m, m')\right)_+ \leq \frac{c}{n}$$

with
$$\tilde{p}(m,m') = 4\frac{D_m \vee D_{m'}}{n} := 4H^2,$$

$$\sup_{t \in B_{m',m}^g(0,1)} \text{Var}[(\delta_1 - F(U_1))t(U_1)] \leq \sup_{t \in B_{m',m}^g(0,1)} \mathbb{E}(t^2(U_1)) = 1 := v,$$

and $\sup_{t \in B_{m',m}^g(0,1)} \|(\delta_1 - F(U_1))t\|_\infty \leq \sup_{t \in B_{m',m}^g(0,1)} \|t\|_\infty \leq \Phi_0\sqrt{D_{m,m'}}/g_0 := M_1.$

## Acknowledgements

Here, we would like to express to Y. Rozenholc our sincere appreciation for providing us very performing Matlab programs. All the material is available on [http://www.math-info.univ-paris5.fr/∼rozen](http://www.math-info.univ-paris5.fr/∼rozen).